\pdfoutput=1
\documentclass{article}

\usepackage{arxiv}

\usepackage{amsmath}
\usepackage{amssymb}
\usepackage{amsthm}
\usepackage{mathtools}
\usepackage{thmtools}
\usepackage{thm-restate}

\usepackage{xcolor}
\usepackage{graphicx}
\usepackage{pifont}
\usepackage{xspace}

\usepackage{algorithmic}
\usepackage{algorithm}

\usepackage{hyperref}
\hypersetup{colorlinks=false}
\usepackage{cleveref}

\definecolor{mydarkgreen}{RGB}{39,130,67}
\definecolor{mydarkred}{RGB}{192,25,25}
\newcommand{\checkmarkgreen}{\textbf{\color{mydarkgreen}\ding{52}}}
\newcommand{\crossmarkred}{{\color{mydarkred}\ding{56}}}

\newcommand{\algname}[1]{{\color{gray}\small\sf#1}\xspace}

\newcommand{\R}{\mathbb{R}}
\newcommand{\cX}{\mathcal{X}}
\newcommand{\norm}[1]{{\left\| #1 \right\|}}
\newcommand{\eqdef}{\coloneqq}
\DeclareMathOperator{\dist}{dist}
\DeclareMathOperator{\dom}{dom}
\newcommand{\broxSS}[3]{{\rm brox}_{#1}^{#2}\left(#3\right)}

\usepackage[utf8]{inputenc}
\usepackage[T1]{fontenc}
\usepackage{booktabs}
\usepackage{microtype}
\usepackage{enumitem}

\usepackage{natbib}

\usepackage{tikz}
\usetikzlibrary{arrows.meta, calc, patterns, shapes.geometric, decorations.pathreplacing}
\usepackage{makecell}

\usepackage{refcount}
\newcommand{\proofof}[2]{Proof of \texorpdfstring{\Cref{#2}}{#1 \getrefnumber{#2}}}

\DeclareMathOperator*{\argsmin}{arg\,min}
\newcommand{\algcomment}[1]{\hfill $\triangleright$ #1}

\theoremstyle{plain}
\newtheorem{theorem}{Theorem}[section]
\newtheorem{assumption}[theorem]{Assumption}

\newtheorem{proposition}[theorem]{Proposition}

\theoremstyle{definition}
\newtheorem{definition}[theorem]{Definition}

\theoremstyle{remark}
\newtheorem{remark}[theorem]{Remark}

\crefname{assumption}{Assumption}{Assumptions}
\Crefname{assumption}{Assumption}{Assumptions}

\crefname{definition}{Definition}{Definitions}
\Crefname{definition}{Definition}{Definitions}

\crefname{remark}{Remark}{Remarks}
\Crefname{remark}{Remark}{Remarks}

\crefname{lemma}{Lemma}{Lemmas}
\Crefname{lemma}{Lemma}{Lemmas}

\crefname{theorem}{Theorem}{Theorems}
\Crefname{theorem}{Theorem}{Theorems}

\crefname{algorithm}{Algorithm}{Algorithms}
\Crefname{algorithm}{Algorithm}{Algorithms}

\title{Convergence Guarantees for the Broximal Point Method under Perturbations}

\author{%
  Alina~Abdikarimova\thanks{Corresponding author: \texttt{alina.abdikarimova@kaust.edu.sa}} \\
  KAUST \\
  \And
  Ammar~Mahran \\
  KAUST \\
  \And
  Peter~Richt\'{a}rik \\
  KAUST \\
}

\begin{document}

\maketitle

\begin{abstract}
The \textit{Broximal Point Method} (\algname{BPM}) of \citet{gruntkowska2025ball} comes with strong guarantees for convex optimization: under convexity alone, it terminates after finitely many steps and converges linearly in function value, with bounds that do not depend on the condition number or on any other regularity property of the objective.
Each iteration, however, requires the exact minimization of the objective over a ball, which is generally intractable.
As a first step toward implementable variants, we study the method when these subproblems are solved only approximately.
We introduce and analyze \algname{Inexact BPM}, which replaces the exact broximal point by an approximation and makes no assumption on how this approximation is obtained.
Under summable or sufficiently small errors, we bound the number of iterations needed to reach a neighborhood of the set of minimizers, and these bounds remain independent of the regularity of the objective.
Under the weaker requirement that the errors stay below the ball radius, we prove sublinear convergence in function value for the exact broximal points underlying the perturbed iterates, valid until the target neighborhood is reached, and we recover linear rates under bounded sublevel sets or strong convexity.
\end{abstract}

\section{Introduction}

We consider the problem
\begin{equation}
\label{eq:problem-def}
    \min_{x\in \R^d} \phi(x),
\end{equation}
where $\phi: \R^d \to \R \cup \{+\infty\}$ is proper, closed, and convex. 
We assume that $\phi$ has a nonempty set of minimizers, that is,
\begin{equation}
    \cX_\star \eqdef \argsmin_{x \in \R^d} \phi(x) \neq \emptyset .
\end{equation}
Since $\cX_\star$ is nonempty, the minimum is attained. We denote the
corresponding optimal value by
\begin{equation}
    \phi_\star \eqdef \min_{x\in\R^d}\phi(x).
\end{equation}
Minimizing a convex function is a fundamental problem in optimization, with applications spanning signal processing, operations research, and the training of large-scale machine learning models. 

A classical algorithmic framework for solving \eqref{eq:problem-def} is the \textit{Proximal Point Method} (\algname{PPM}), originally introduced by \citet{martinet1970regularisation} and later studied by \citet{rockafellar1976monotone}. 
\algname{PPM} is an iterative scheme in which a quadratic penalty term is added to the convex objective function. 
Starting from a point $x_k$, the next iterate $x_{k+1}$ is the solution of the following minimization problem,
\begin{equation}
\label{eq:ppm-step}
    x_{k+1} = \argsmin_{x \in \R^d} \left\{ \phi(x) + \frac{1}{2\gamma_k}\norm{x-x_k}^2 \right\} ,
\end{equation}
where $\gamma_k > 0$ acts as a stepsize parameter that controls the trade-off between minimizing the objective function and remaining close to the current iterate.
By adding the quadratic penalty term to the objective, \algname{PPM} requires a strongly convex minimization subproblem to be solved at every iteration. 
For nonsmooth convex functions, \algname{PPM} is known to converge sublinearly \citep{guler1991convergence}, whereas linear convergence is achieved under the additional assumption of strong convexity \citep{rockafellar1976monotone}.

While \algname{PPM} relies on a soft quadratic regularization, the \textit{Broximal Point Method} (\algname{BPM}), recently introduced by \citet{gruntkowska2025ball}, utilizes a hard ball constraint. 
The next iterate of the \algname{Exact BPM} algorithm (cf. \Cref{alg:exact-brox}) is the exact minimizer of the objective function over the ball $B_t(x_k)$ of radius $t$ centered at the current iterate $x_k$,
\begin{equation}
\label{eq:bpm-step}
    x_{k+1}\in\argsmin_{x\in B_t(x_k)}\phi(x).
\end{equation}
Under additional differentiability and smoothness assumptions, this update can be viewed as normalized gradient descent on the ball envelope~\citep{gruntkowska2025ball}.
The hard constraint yields stronger guarantees than those available for \algname{PPM}.
Under convexity alone, \algname{Exact BPM} reaches a minimizer after finitely many steps, and its iteration bound depends only on the initial distance to the set of minimizers, not on the condition number or any other regularity property of the objective.
Its function values also converge linearly, and the contraction factor improves as the iterates approach the set of minimizers \citep{gruntkowska2025ball, richtarik2026sharper}.

The broximal framework connects to trust-region methods~\citep{conn2000trust}, normalized gradient descent~\citep{nesterov1984minimization}, and ball-oracle complexity~\citep{carmon2020acceleration, carmon2021thinking}.
In particular, \citet[Appendix~A]{carmon2020acceleration} analyze the iteration of an exact ball-optimization oracle, which is \algname{Exact BPM} in our notation, for strictly convex objectives.
\citet{gruntkowska2025non} extend the convergence results to non-Euclidean geometries.
This makes the framework relevant to recent deep learning optimizers with a similar non-Euclidean structure, such as \algname{Muon}~\citep{jordan2024muon}, \algname{Scion}~\citep{pethick2025training}, and \algname{Gluon}~\citep{riabinin2026muon}, all of which can be interpreted as approximate broximal updates applied to surrogate subproblems~\citep{gruntkowska2025non}.

Neither \algname{BPM} nor \algname{PPM} is directly implementable because the proximal and ball-constrained subproblems \eqref{eq:ppm-step} and \eqref{eq:bpm-step} rarely admit a closed-form solution. 
In the case of \algname{PPM}, this has not hindered the development of implementable algorithms that substitute the exact proximal step with an approximation, such as proximal-gradient methods \citep{parikh2014proximal} and operator splitting techniques \citep{eckstein1992douglas, boyd2011distributed}.

These implementable methods rely on the underlying framework's robustness to computational errors. 
Early work by \citet{rockafellar1976monotone} showed that global convergence of \algname{PPM} is preserved when the approximation errors form a summable sequence. 
This was extended to more general error models \citep{guler1992new, solodov1999hybrid} and stochastic settings \citep{ryu2017stochastic, bianchi2016ergodic, richtarik2024unified, tovmasyan2026revisiting}.
For implementable algorithms, \citet{schmidt2011convergence} established that the convergence rates of proximal-gradient methods are preserved under controlled deterministic errors in the gradient or proximal evaluations.

Motivated by these extensions of \algname{PPM}, we study how \algname{BPM} behaves under similar perturbations.
We introduce \algname{Inexact BPM} (cf. \Cref{alg:inexact-brox}), in which each subproblem is solved approximately, so that the iterates deviate from the exact broximal points.
We impose no structure on how the approximation is obtained, so our guarantees hold for any inner solver that returns a point within the prescribed distance of the exact broximal point.
We leave the design of solver-specific implementations to future work.

Our main question is which guarantees of \algname{Exact BPM} survive such perturbations, in particular finite termination and complexity bounds that do not depend on the regularity of the objective.
Such guarantees are not available for classical \algname{PPM} under convexity alone.
We show that finite, regularity-independent hitting-time bounds persist under summable or sufficiently small errors.
Under persistent errors bounded by a fraction of the ball radius ($\bar{\varepsilon}<t$), we obtain a sublinear function-value bound along the sequence of exact broximal points, valid until the iterates enter the target neighborhood, with an exponent that depends only on the ratio $\bar{\varepsilon}/t$.
This regime is covered neither by classical summable-error analyses of \algname{PPM} nor by the high-accuracy accelerated guarantee of \citet{carmon2020acceleration}.

In \Cref{sec:small-errors}, we show that \algname{Inexact BPM} enters the $t$-tube around the set of minimizers after a bounded number of iterations when the errors are summable or uniformly bounded by a constant that depends on the initial distance to the set of minimizers.
In \Cref{sec:large-errors}, we consider a relaxed uniform bound on the errors and establish sublinear convergence in function value, at a rate that improves as the error bound tightens.
In the limit of vanishing errors, we recover a linear rate for \algname{Exact BPM}.
We also show that linear convergence is recovered under bounded sublevel sets or strong convexity.
\Cref{sec:num-exp} illustrates the results on a simple example.
Proofs not given in the main text are deferred to \Cref{sec:proofs}.

\section{Preliminaries}
\label{sec:prelim-alg}

We work under the following assumption on the objective function $\phi$.
\begin{assumption}
\label{assumption:main}
The function $\phi: \R^d \to \R \cup \{+\infty\}$ is proper, closed, and convex, and its set of minimizers,
\begin{equation}
    \label{eq:def-minimizer-set}
    \cX_\star = \argsmin_{x \in \R^d} \phi(x),    
\end{equation}
is nonempty. 
\end{assumption}

Closedness and convexity of the objective function $\phi$ imply closedness and convexity of its set of minimizers $\cX_\star$, ensuring uniqueness of projections onto it.
We denote the Euclidean distance of an iterate $x_k \in \R^d$ to this set by
\begin{equation}
    D_k \eqdef \dist(x_k, \cX_\star) .
\end{equation}

\subsection{Exact broximal point method}
\label{subsec:exact_brox}

The broximal operator of \citet{gruntkowska2025ball} returns the solutions of the subproblem \eqref{eq:bpm-step}.

\begin{definition}[Broximal operator]
\label{def:brox-operator}
For a radius $t>0$ and a function $\phi$ satisfying \Cref{assumption:main}, the broximal operator is the set-valued mapping
\begin{equation}
    \label{eq:def-brox-operator}
    \broxSS{\phi}{t}{x} \eqdef \argsmin_{y \in B_t(x)} \phi(y),
\end{equation}
where $B_t(x) \eqdef \{y \in \R^d : \norm{y - x} \le t\}$ is the closed Euclidean ball of radius $t$ centered at $x$.
\end{definition}

Evaluating the broximal operator thus requires solving a constrained optimization problem.
While this operator is set-valued, the resulting set is typically a singleton.

\begin{remark}[Properties of the broximal operator]
    \label{rem:brox-properties}
    Under \Cref{assumption:main}, if $B_t(x)\cap\dom\phi\neq\emptyset$, then $\broxSS{\phi}{t}{x}$ is nonempty and the subproblem has a finite minimum. 
    Moreover, if $B_t(x) \cap \cX_\star = \emptyset$, that is, if $\dist(x, \cX_\star) > t$, then $\broxSS{\phi}{t}{x}$ is a singleton $\{y\}$ lying on the boundary of the ball:
    \begin{equation}
        \label{eq:brox-boundary}
        \norm{y - x} = t.
    \end{equation}

    In such cases, we write $y = \broxSS{\phi}{t}{x}$ for this unique point.
    These properties are shown in the proof of \Cref{lem:exact-one-step}; see also \citet[Theorem~D.1]{gruntkowska2025ball}.
\end{remark}

The following lemma establishes that the broximal operator decreases the distance to the set of minimizers.

\begin{restatable}[One-step contraction {\cite[Theorem~8.1(ii)]{gruntkowska2025ball}}]{lemma}{broxlemma}
    \label{lem:exact-one-step}
    Let \Cref{assumption:main} hold. 
    For any $x \in \R^d$ satisfying $B_t(x)\cap\dom\phi\neq\emptyset$ and $\dist(x, \cX_\star) > t$, let $y = \broxSS{\phi}{t}{x}$ be the unique broximal point (cf. \Cref{rem:brox-properties}). 
    Then, for every minimizer $x_\star \in \cX_\star$, we have:
    \begin{equation}
        \label{eq:exact-one-step}
        \norm{y - x_\star}^2 \le \norm{x - x_\star}^2 - t^2.
    \end{equation}
\end{restatable}

\begin{algorithm}[!t]
    \begin{algorithmic}[1]
        \STATE \textbf{Input:} Starting point $x_0 \in \dom \phi$, radius $t > 0$
        \FOR{$k=0,1,2,\ldots$}
            \IF{$\dist(x_k,\cX_\star) > 0$}
                \STATE $x_{k+1}\in\broxSS{\phi}{t}{x_k}$
            \ELSE
                \STATE $x_{k+1} = x_k$
            \ENDIF
        \ENDFOR
    \end{algorithmic}
    \caption{\algname{Exact BPM}: The exact broximal point method.}
    \label{alg:exact-brox}
\end{algorithm}

By \Cref{lem:exact-one-step}, as long as $D_k > t$, each exact broximal step reduces the squared distance to $\cX_\star$ by at least $t^2$; once a minimizer lies within distance $t$ of $x_k$, the next step reaches it.
This argument uses only \Cref{assumption:main} and requires no differentiability, smoothness, finite-valuedness, or strong convexity.
Consequently, starting from $x_0 \in \dom \phi$, \algname{Exact BPM} reaches a minimizer in at most
\begin{equation}
    \label{eq:iteration-count-exact-brox}
    \left\lceil \frac{D_0^2}{t^2} \right\rceil
\end{equation}
iterations \citep[Theorem~8.1(iii)]{gruntkowska2025ball}, where $D_0 = \dist(x_0, \cX_\star)$ is the distance of the initial point to the set of minimizers.
Classical methods for general convex optimization, such as \algname{PPM} \citep{guler1991convergence}, subgradient methods \citep{shor1985minimization, nesterov2018lectures}, and proximal-gradient methods \citep{parikh2014proximal}, come with convergence rates but do not in general terminate after finitely many steps; the bound \eqref{eq:iteration-count-exact-brox} is finite and depends only on $D_0$ and $t$, not on smoothness, the condition number, or any other structural parameter of the objective.
As the dimension grows, this dimension-independent bound is asymptotically attained by a family of polyhedral functions \citep[Theorem~4.11]{richtarik2026sharper}; in fixed dimension, the worst-case iteration count grows only linearly in $D_0/t$ \citep[Theorem~4.5]{richtarik2026sharper}.

\subsection{Inexact broximal point method}

Each iteration of \algname{Exact BPM} (\Cref{alg:exact-brox}) requires minimizing the objective function over a ball, which is generally intractable in practice. 
We therefore consider an inexact variant where the exact broximal point is replaced by an approximation.
Specifically, we assume that the next iterate $x_{k+1}$ produced by \algname{Inexact BPM} (\Cref{alg:inexact-brox}) lies no more than some distance $\varepsilon_k \ge 0$ from the exact broximal step evaluated at the current iterate $x_k$, that is,
\begin{equation}
    \label{eq:approximation-error}
    \norm{x_{k+1}-\broxSS{\phi}{t}{x_k}}\le\varepsilon_k,
    \qquad \text{whenever }D_k>t.
\end{equation}
We impose no structural assumptions on how this approximation is obtained,
and our framework thus accommodates deterministic, randomized, and even adversarial perturbations.
We consider a variety of settings for the error tolerances $\{\varepsilon_k\}$ below.

\textbf{Relation to inexact ball oracles.}
The error model \eqref{eq:approximation-error} measures inexactness by the distance from the produced iterate to the exact ball-constrained minimizer.
\citet{carmon2020acceleration} use the same notion in their $(\delta, r)$-ball optimization oracle, which, when queried at $\bar{x}$, returns a point $y \in B_r(\bar{x})$ with $\norm{y - x_{\bar{x}, r}} \le \delta$ for some $x_{\bar{x}, r} \in \argsmin_{x \in B_r(\bar{x})} \phi(x)$.
Their setting differs from ours in several ways.
Their accelerated guarantee assumes real-valued, strictly convex, smooth objectives and uses the oracle inside Monteiro--Svaiter acceleration, with an accuracy $\delta$ that must be a small fraction of $r$, depending on the smoothness constant, the initial distance to the optimum, and the target accuracy.
We require only that $\phi$ be proper, closed, and convex, we analyze the broximal iteration itself, and we allow persistent errors as large as a constant fraction of the radius ($\bar{\varepsilon} < t$; cf.\ \Cref{assumption:relaxed-error}).
Moreover, their oracle output lies in $B_r(\bar{x})$, whereas our iterate $x_{k+1}$ need not lie in $B_t(x_k)$.

With approximation errors, the iterates can in general only be guaranteed to reach a neighborhood of $\cX_\star$, not $\cX_\star$ itself.
We use the following neighborhood as the target.
\begin{definition}[Tube]
    \label{def:t-tube}
    Given $\rho \ge 0$, the closed $\rho$-tube around a set $\cX \subseteq \R^d$ is
    \begin{equation}
        \label{eq:t-tube-def}
        \mathcal{T}_\rho(\cX) \eqdef \bigl\{ x \in \R^d : \dist(x, \cX) \le \rho \bigr\}.
    \end{equation}
\end{definition}

Since $\cX_\star$ is closed, $\mathcal{T}_0(\cX_\star) = \cX_\star$, so the $t$-tube $\mathcal{T}_t(\cX_\star)$ generalizes the target set of the exact method.
The radius $t$ thus plays two roles: it is the radius of the ball constraint and the target accuracy.

\Cref{alg:inexact-brox} generalizes \Cref{alg:exact-brox} to inexact steps.
It includes a truncation step that freezes the sequence once the iterates enter $\mathcal{T}_t(\cX_\star)$.
This truncation is an analytical device, not an implementable stopping rule: the test $\dist(x_k, \cX_\star) \le t$ cannot be evaluated in general, since $\cX_\star$ is unknown.
The updates themselves do not require the test.

\begin{algorithm}[ht]
    \begin{algorithmic}[1]
        \STATE \textbf{Input:} Starting point $x_0 \in \dom \phi$, radius $t>0$
        \FOR{$k=0,1,2,\ldots$}
            \IF{$\dist(x_k,\cX_\star) > t$}
                \STATE $x_{k+1} \approx \broxSS{\phi}{t}{x_k}$ 
                \algcomment{The approximation satisfies \eqref{eq:approximation-error}.}
            \ELSE
                \STATE $x_{k+1} = x_k$
                \algcomment{Truncation step.}
            \ENDIF
        \ENDFOR
    \end{algorithmic}
    \caption{\algname{Inexact BPM}: The truncated inexact broximal point method.}
    \label{alg:inexact-brox}
\end{algorithm}

Compared with \Cref{alg:exact-brox}, \Cref{alg:inexact-brox} makes two changes.
First, $x_{k+1}$ only approximates $\broxSS{\phi}{t}{x_k}$, with error at most $\varepsilon_k \ge 0$ as in \eqref{eq:approximation-error}; the tolerances need not be known to run the method.
Second, the method stops when it enters $\mathcal{T}_t(\cX_\star)$ rather than when it reaches $\cX_\star$.

To quantify the convergence properties of \algname{Inexact BPM}, we analyze the number of iterations required for the sequence to enter the $t$-tube $\mathcal{T}_t(\cX_\star)$. 
We formalize this through the notion of a hitting time.

\begin{definition}[Hitting time]
    \label{def:hitting-time}
    The hitting time $K_t$ of \Cref{alg:inexact-brox} to the $t$-tube $\mathcal{T}_t(\cX_\star)$ is defined as:
    \begin{equation}
    \label{eq:hitting-time-def}
        K_t \eqdef \inf \bigl\{ k \ge 0 : x_k \in \mathcal{T}_t(\cX_\star) \bigr\} = \inf \bigl\{ k \ge 0 : \dist(x_k, \cX_\star) \le t \bigr\},
    \end{equation}
    where we adopt the standard convention that $\inf \emptyset = +\infty$.
\end{definition}

Without errors, $K_t$ is bounded by the iteration count \eqref{eq:iteration-count-exact-brox} of the exact method.
In the next section, we derive analogous bounds under perturbations.
When the approximations are randomized, $K_t$ is a random variable; 
all our bounds are worst-case and hold pathwise, for every sequence of approximations that satisfies the stated error conditions.

\section{Distance contraction guarantees}
\label{sec:small-errors}

In this section, we establish the iteration complexity of \algname{Inexact BPM} under well-controlled error regimes. 
To that end, we first derive a one-step distance recursion generalizing \Cref{lem:exact-one-step} for the inexact setting. 

\begin{restatable}{lemma}{recursionlemma}
\label{lem:dist-recursion}
    Let $\{x_k\}_{k\ge 0}$ be the sequence generated by \Cref{alg:inexact-brox}.
    Suppose \Cref{assumption:main} holds and $B_t(x_k)\cap\dom\phi\neq\emptyset$ for every $k<K_t$. 
    The distances $D_k = \dist(x_k, \cX_\star)$ to the set of minimizers satisfy the following relations:
    \begin{enumerate}[label=(\roman*)]
        \item
        If $D_k > t$, we have the recursive upper bound:
        \begin{equation}
        \label{eq:basic_recursion}
            D_{k+1} \le \sqrt{D_k^2 - t^2} + \varepsilon_k.
        \end{equation}
        
        \item
        If $D_k \le t$, the algorithm truncates the update step ($x_{k+1} = x_k$), and thus:
        \begin{equation}
        \label{eq:truncation_recursion}
            D_{k+1} = D_k \le t.
        \end{equation}
    \end{enumerate}
\end{restatable}

Thus the contraction of the exact step in \Cref{lem:exact-one-step} survives, up to an additive error $\varepsilon_k$ in the distance to the set of minimizers.
We now consider two classical error regimes in which this recursion yields finite hitting times.

\subsection{Summable errors}

We first allow individual errors to be large, provided their sum is finite.
This is the error model of the classical inexact proximal point analysis of \citet{rockafellar1976monotone}.

\begin{restatable}{theorem}{summable}
\label{theorem:summable-progress-guarantee}
    Let $\{x_k\}_{k\ge 0}$ be the sequence generated by \Cref{alg:inexact-brox} with radius $t > 0$, and suppose \Cref{assumption:main} holds.
    Assume further that every active subproblem is well-defined, in the sense that $B_t(x_k) \cap \dom \phi \neq \emptyset$ for all $k < K_t$.
    If the error sequence is summable, meaning
    \begin{equation}
    \label{eq:summable-error-def}
        \Sigma \eqdef \sum_{k=0}^\infty \varepsilon_k < \infty,
    \end{equation}
    then the hitting time $K_t$ to the $t$-tube $\mathcal{T}_t(\cX_\star)$ satisfies:
    \begin{equation}
    \label{eq:summable-hitting-time-bound}
        K_t \le \left\lceil \frac{(D_0 + \Sigma)^2}{t^2} \right\rceil.
    \end{equation}
\end{restatable}

When $\Sigma=0$, the bound reduces to $\lceil D_0^2/t^2\rceil$, the same expression as the \algname{Exact BPM} baseline \eqref{eq:iteration-count-exact-brox}, although $K_t$ measures tube entry rather than arrival at a minimizer.
The term $(D_0+\Sigma)^2$ accounts for the worst-case accumulation of approximation errors, which may increase the distance to the set of minimizers by up to $\Sigma$ relative to its initial value $D_0$.
Thus finite tube entry is preserved under summable errors, with the worst-case hitting-time bound growing quadratically in the total error budget $\Sigma$.

If, for example, $\Sigma\le cD_0$ for some $c\ge0$, then \eqref{eq:summable-hitting-time-bound} gives $K_t\le\lceil(1+c)^2D_0^2/t^2\rceil$.

The condition $B_t(x_k) \cap \dom \phi \neq \emptyset$ in \Cref{theorem:summable-progress-guarantee} is needed because the summable regime places no bound on any individual error, so a single large perturbation may move $x_k$ to a point whose ball misses $\dom\phi$; the condition holds automatically when $\varepsilon_k < t$ for all $k$ (see the proof of \Cref{theorem:uniform-small-errors}).

\subsection{Uniformly small errors}

We next consider errors that are uniformly bounded across all iterations.
To quantify the progress made at each step, we use the following quantity.

\begin{definition}[Progress guarantee]
\label{def:progress-guarantee}
    Let $\{x_k\}_{k \ge 0}$ be the sequence generated by \Cref{alg:inexact-brox} with radius $t$, and let $\{\varepsilon_k\}_{k \ge 0}$ denote a sequence of nonnegative error bounds for the approximations as described in \eqref{eq:approximation-error}. 
    For each iteration $k$ satisfying $D_k > t$, the progress guarantee $p_k$ is defined as:
    \begin{equation}
    \label{eq:progress-guarantee-def}
        p_k \eqdef t^2 - \varepsilon_k^2 - 2 \varepsilon_k \sqrt{D_k^2 - t^2}.
    \end{equation}
\end{definition}

Since $\varepsilon_k \ge 0$, we have $p_k \le t^2$.
Squaring both sides of the distance recursion in \Cref{lem:dist-recursion}(i) yields
\begin{equation}
\label{eq:progress-squared-relation}
    D_{k+1}^2 \le D_k^2 - p_k,
\end{equation}
so the progress guarantee $p_k$ quantifies the minimum guaranteed reduction in the squared distance to the set of minimizers $\cX_\star$ at step $k$. 
For $\varepsilon_k = 0$, we recover $p_k = t^2$, the exact per-step reduction of \algname{Exact BPM}.

The progress guarantee $p_k$ depends on the radius $t$, the error bound $\varepsilon_k$, and the current distance $D_k$.
It is positive only if the error bound is small enough.
Specifically, we have $p_k > 0$ if and only if
\begin{equation}
\label{eq:error-positivity-threshold}
    \varepsilon_k < D_k - \sqrt{D_k^2 - t^2}.
\end{equation}
If this condition is violated ($p_k \le 0$), the approximation error may exceed the broximal contraction, potentially pushing the iterate further away from the set of minimizers. 
Under a uniform bound on the errors that prevents the iterates from moving away from the set of minimizers, we obtain the following guarantee for the hitting time.

\begin{restatable}{theorem}{uniformsmall}
\label{theorem:uniform-small-errors}
    Let $\{x_k\}_{k\ge 0}$ be the sequence generated by \Cref{alg:inexact-brox} with radius $t < D_0$, and suppose \Cref{assumption:main} holds. 
    Suppose the error sequence is uniformly bounded such that $\varepsilon_k \le \bar{\varepsilon}$ for all $k \ge 0$, and let
    \begin{equation}
    \label{eq:p0-definition}
        \bar{p} \eqdef t^2 - \bar{\varepsilon}^2 - 2\bar{\varepsilon}\sqrt{D_0^2 - t^2}
    \end{equation}
    denote the initial worst-case progress guarantee. 
    If the uniform error bound ensures that $\bar{p} > 0$, that is, if $\bar{\varepsilon} < D_0 - \sqrt{D_0^2 - t^2}$, then the hitting time $K_t$ to the $t$-tube $\mathcal{T}_t(\cX_\star)$ satisfies:
    \begin{equation}
    \label{eq:uniform-hitting-time-bound-p0}
        K_t \le \left\lceil \frac{D_0^2}{\bar{p}} \right\rceil.
    \end{equation}
\end{restatable}

\begin{table}[t]
\centering
\caption{\textbf{Hitting time bounds.}
The hitting time $K_t$ to the $t$-tube is defined in \eqref{eq:hitting-time-def}.
Here $D_0 = \dist(x_0, \cX_\star)$ is the initial distance to the set of minimizers, and $\bar{p} = t^2 - \bar{\varepsilon}^2 - 2\bar{\varepsilon}\sqrt{D_0^2 - t^2} \le t^2$ is the initial worst-case progress guarantee \eqref{eq:p0-definition} under a uniform error bound $\bar{\varepsilon}$.
All bounds are independent of the condition number and of any other regularity parameter of the objective.}
\label{tab:sec3}
\small
\begin{tabular}{ccc}
\toprule
\textbf{Result} & \textbf{Error regime} & \textbf{Bound on} $K_t$ \\
\midrule
\makecell{Theorem 8.1(iii) \\ \citep{gruntkowska2025ball}} & \makecell{No errors: \\ $\varepsilon_k = 0$} & $\left\lceil \frac{D_0^2}{t^2} \right\rceil$ \\
\addlinespace
\Cref{theorem:summable-progress-guarantee} & \makecell{Summable errors: \\ $\Sigma = \sum_{k=0}^\infty \varepsilon_k < \infty$} & $\left\lceil \frac{(D_0 + \Sigma)^2}{t^2} \right\rceil$ \\
\addlinespace
\Cref{theorem:uniform-small-errors} & \makecell{Uniformly small errors: \\ $\varepsilon_k \leq \bar{\varepsilon} < D_0 - \sqrt{D_0^2 - t^2}$} & $\left\lceil \frac{D_0^2}{\bar{p}} \right\rceil$ \\
\bottomrule
\end{tabular}
\end{table}

The bound \eqref{eq:uniform-hitting-time-bound-p0} has the same form as the exact baseline \eqref{eq:iteration-count-exact-brox}, with the per-step progress $t^2$ replaced by $\bar{p}$.
The proof shows that $p_k \ge \bar{p}$ for all $k < K_t$, however the approximation errors evolve, so the squared distance to $\cX_\star$ decreases by at least $\bar{p}$ in every iteration before tube entry.
Like \eqref{eq:iteration-count-exact-brox}, the bound does not depend on any property of the objective.
\Cref{tab:sec3} summarizes the results of this section.

The threshold in \eqref{eq:error-positivity-threshold} decreases in $D_k$, and for any fixed radius $t$,
\begin{equation}
    \lim_{D_0 \to +\infty} \Bigl(D_0 - \sqrt{D_0^2 - t^2}\Bigr) = 0 .
\end{equation}
That is, for any fixed $\bar{\varepsilon} > 0$, a sufficiently distant starting point gives $\bar{p} < 0$, and no single uniform tolerance guarantees positive progress from every starting point.
This limitation motivates the next section.

\section{Suboptimality bounds under relaxed tolerances}
\label{sec:large-errors}

We now only require the errors to stay strictly below the ball radius.

\begin{assumption}[Radius-bounded uniform error]
    \label{assumption:relaxed-error}
    The error bounds $\{\varepsilon_k\}_{k \ge 0}$ in \eqref{eq:approximation-error} are uniformly bounded by a constant $\bar{\varepsilon}$ smaller than the radius $t$:
    \begin{equation}
    \label{eq:relaxed-error-bound}
        \varepsilon_k \le \bar{\varepsilon} < t, \quad \forall k \ge 0.
    \end{equation}
\end{assumption}

The error bound $\bar{\varepsilon}$ may be arbitrarily close to $t$.
This condition is weaker than the one in \Cref{theorem:uniform-small-errors} and does not involve the distance to the set of minimizers.
Unlike the summable model, which forces the errors to vanish, \Cref{assumption:relaxed-error} allows a persistent error at every iteration.

The strict inequality $\bar{\varepsilon} < t$ cannot be dropped.
Since $\norm{\broxSS{\phi}{t}{x_k} - x_k} \le t$, an adversary with $\bar{\varepsilon} \ge t$ can choose $x_{k+1} = x_k$ at every step, so the method stalls at $x_0$.

Throughout this section, $t$ and $\bar{\varepsilon}$ are fixed constants satisfying \Cref{assumption:relaxed-error}.
We define the inner radius $r$,
\begin{equation}
\label{eq:r-and-R}
    r \eqdef t - \bar{\varepsilon} ,
\end{equation}
and the auxiliary sequence $\{y_k\}$,
\begin{equation}
    \label{eq:y-seq}
    \begin{aligned}
        y_0 &\eqdef x_0,\\
        y_{k+1} &\eqdef \broxSS{\phi}{t}{x_k}, \qquad 0\le k<K_t.
    \end{aligned}
\end{equation}
The point $y_{k+1}$ is the point that \algname{Inexact BPM} would compute in iteration $k$ without approximation error.
With this notation, the approximation error bound \eqref{eq:approximation-error} may equivalently be expressed as
\begin{equation}
    \label{eq:approximation-error-y}
    \norm{x_{k+1}-y_{k+1}}\le\varepsilon_k,
    \qquad 0\le k<K_t.
\end{equation}
Equivalently, for every finite $1\le k\le K_t$,
\begin{equation}
\label{eq:approximation-error-general}
    \norm{x_k-y_k}
    \le \varepsilon_{k-1}
    \le \bar{\varepsilon}.
\end{equation}
Since $y_0=x_0$, the bound $\norm{x_k-y_k}\le\bar{\varepsilon}$ holds for every finite $0\le k\le K_t$.

\begin{figure}[htbp]
\centering
\begin{tikzpicture}[scale=1.1, >=stealth, even odd rule]

    \coordinate (yk) at (0,0);
    \coordinate (xstar) at (4.5, 2.0);
    \coordinate (xk) at (-0.2, -0.4); 
    
    \def\rt{2.5}
    \def\reps{0.6}
    \def\rr{1.9}
    
    \draw[thick, fill=gray, fill opacity=0.08, draw=gray] (xk) circle (\rt);
    \node[gray] at (-2.7, 1.5) {$B_t(x_k)$};
    
    \draw[thick, dashed, fill=blue, fill opacity=0.05, draw=blue] (yk) circle (\rr);
    \node[blue] at (0.9, -1.0) {$B_r(y_k)$};
    
    \draw[fill=red, fill opacity=0.12, draw=red] (yk) circle (\reps);
    \node[red, scale=0.85] at (-1.0, 0.4) {$B_{\bar{\varepsilon}}(y_k)$};
    
    \draw[ultra thick, gray] ($(xstar)+(-0.4,0.8)$) to[out=-60, in=120] ($(xstar)+(0.4,-0.8)$);
    \node[gray, right=6pt, yshift=4pt] at (xstar) {$\cX_\star$};

    \draw[line width=0.8pt, gray, dashed] (yk) -- (xstar);
    
    \coordinate (zk) at ($(yk)!{\rr/sqrt(4.5^2 + 2.0^2)}!(xstar)$);
    
    \coordinate (yk1) at ($(xk) + (15:\rt)$);

    \filldraw[black] (xk) circle (1.5pt) node[below left] {$x_k$};
    \filldraw[black] (yk) circle (1.5pt) node[above=2pt, scale=0.95] {$y_k$};
    \filldraw[black] (xstar) circle (1.5pt) node[above=2pt] {};
    \node[black] at (4.15, 2.15) {$y_k^\star$};
    \filldraw[blue] (zk) circle (1.8pt) node[above left=2pt] {$z_k$};
    \filldraw[black] (yk1) circle (1.8pt) node[right=3pt, scale=0.95] {$y_{k+1}$};

    \draw[->, gray] (xk) -- +(225:\rt) node[midway, left, scale=0.9] {};
    \draw[black,
        decorate,
        decoration={brace, amplitude=3pt, raise=4pt}]
    ($(xk)+(225:0.05)$) -- ($(xk)+(225:\rt-0.05)$);
    \node at ($(xk)+(225:\rt/2)+(0.27,-0.25)$) {$t$};
    \draw[->, blue, thick] (yk) -- (zk) node[midway, above left=-1pt, scale=0.9] {};
    \draw[black,
        decorate,
        decoration={brace, mirror, amplitude=3pt, raise=4pt}]
    (yk) -- (zk);
    \node at ($(yk)!0.5!(zk) + (0.2,-0.35)$) {$r$};
    \draw[->, red, thick] (yk) -- (xk) node[midway, left=1pt, scale=0.85] {};

\end{tikzpicture}
  \caption{\textbf{Geometry of \Cref{assumption:relaxed-error}.} The approximation error bound $\norm{x_k - y_k} \le \bar{\varepsilon}$ guarantees the inclusion $B_r(y_k) \subseteq B_t(x_k)$ for $r = t - \bar{\varepsilon}$. 
  Thus, a reference point $z_k$ on the segment $[y_k, y_k^\star]$, where $y_k^\star$ is the projection of $y_k$ onto $\cX_\star$, lies in $B_t(x_k)$ if $\norm{y_k - z_k} \le r$.
  Since $y_{k+1}$ minimizes $\phi$ over $B_t(x_k)$, we get $\phi(y_{k+1}) \le \phi(z_k)$, which yields a decrease in function value.}
  \label{fig:geometry}
\end{figure}
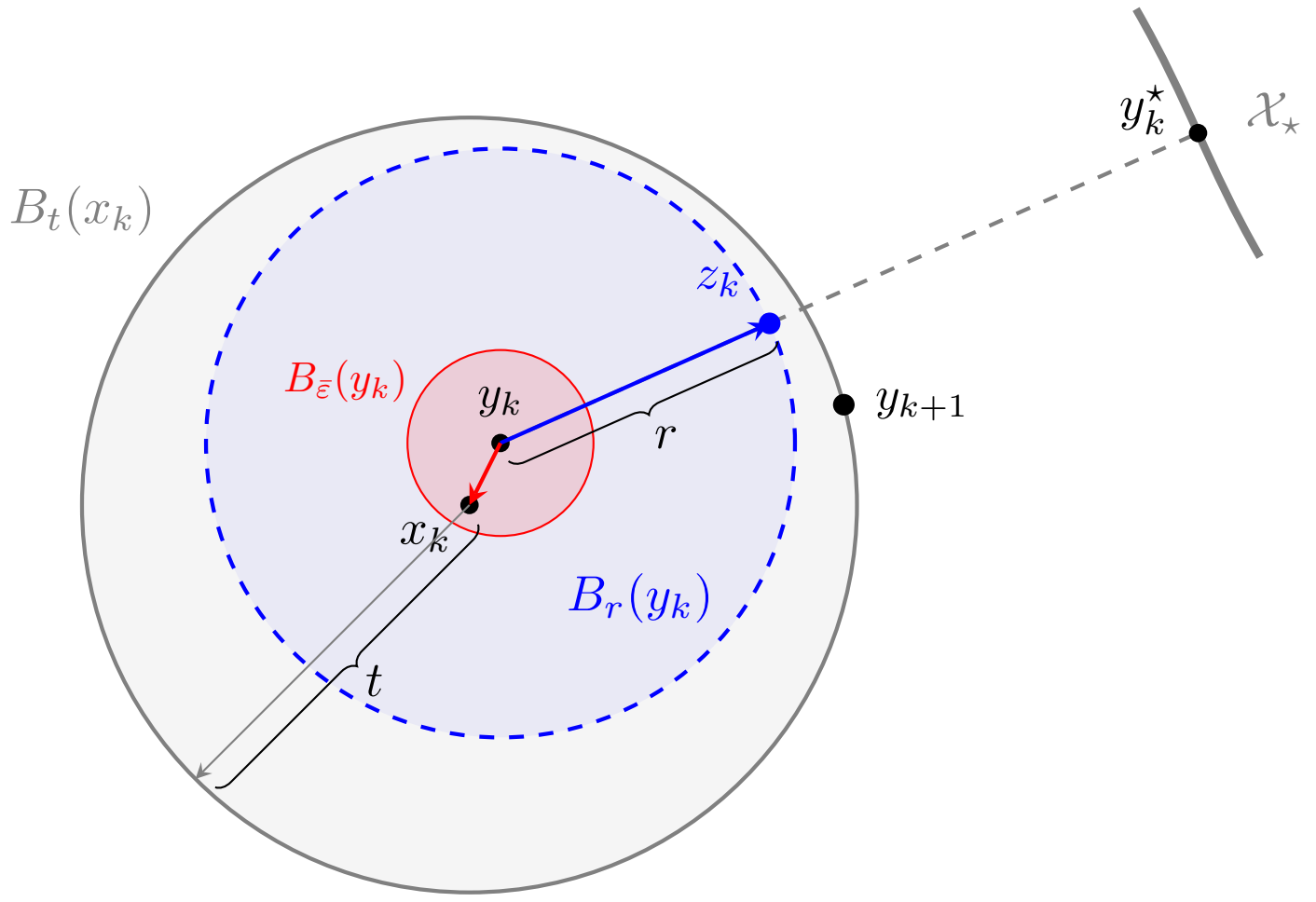

For general convex objectives, we have no bound on $K_t$ under \Cref{assumption:relaxed-error} alone, but it is finite under bounded sublevel sets (see the discussion after \Cref{theorem:bounded-sublevel}).
Moreover, the distance-based condition \eqref{eq:approximation-error-y} does not control the objective value $\phi(x_k)$ under \Cref{assumption:main} alone.
We therefore analyze the decrease of the objective along the sequence $\{y_k\}$ of exact broximal points.
The key observation, illustrated in \Cref{fig:geometry}, is that $\norm{x_k - y_k} \le \bar{\varepsilon}$ implies $B_r(y_k) \subseteq B_t(x_k)$: any point within distance $r$ of $y_k$ lies within distance $t$ of $x_k$ and is therefore feasible for the subproblem that defines $y_{k+1} = \broxSS{\phi}{t}{x_k}$.

\subsection{General convex objectives}

Our first result requires no assumption on $\phi$ beyond \Cref{assumption:main}.
It shows that the suboptimality gap decays at a sublinear rate governed by the ratio $\bar{\varepsilon}/t$.

\begin{restatable}{theorem}{generalconvex}
\label{theorem:two-step-analysis}
    Let $\{x_k\}_{k\ge 0}$ be the sequence generated by \Cref{alg:inexact-brox} and $\{y_k\}_{k\ge 0}$ the auxiliary sequence defined in \eqref{eq:y-seq}, and suppose \Cref{assumption:main,assumption:relaxed-error} hold. 
    Assume $D_0 > t$ and that the uniform error bound of \Cref{assumption:relaxed-error} is positive, $\bar{\varepsilon} > 0$. 
    For any iteration $N \le K_t$, the suboptimality gap satisfies:
    \begin{equation}
        \label{eq:general-convex-rate}
        \phi(y_N) - \phi_\star \le \left(1 + \frac{N \bar{\varepsilon}}{D_0}\right)^{-\frac{r}{\bar{\varepsilon}}}(\phi(y_0) - \phi_\star),
    \end{equation}
    where $r$ is the inner radius defined in \eqref{eq:r-and-R}.
\end{restatable}

\begin{proof}
    The following proof adapts the auxiliary point construction of \citet[Theorem~22]{carmon2020acceleration}.
    Fix an iteration $k < K_t$.
    Since $k<K_t$, we have $\dist(x_k,\cX_\star)>t$. Therefore,
    \[
        \dist(y_k,\cX_\star) \ge \dist(x_k,\cX_\star)-\norm{x_k-y_k} > t-\bar{\varepsilon} = r.
    \]
    Let $y_k^\star$ be the unique projection of $y_k$ onto
    $\cX_\star$. Then
    \[
        \norm{y_k^\star-y_k} = \dist(y_k,\cX_\star) > r,
    \]
    and hence
    \[
        \alpha \eqdef\frac{r}{\norm{y_k^\star-y_k}} \in(0,1).
    \]
    Define the reference point
    \begin{equation}
        z_k \eqdef y_k + \frac{r}{\norm{y_k^\star - y_k}} (y_k^\star - y_k).
    \end{equation}
    Then $z_k = (1-\alpha) y_k + \alpha y_k^\star$ lies at distance $r$ from $y_k$ on the line segment connecting $y_k$ and $y_k^\star$.
    Because $\phi(y_k^\star) = \phi_\star$, convexity of $\phi$ yields
    \begin{equation}
        \label{eq:general-convexity}
        \phi(z_k)-\phi_\star \le (1-\alpha)\bigl(\phi(y_k)-\phi_\star\bigr).
    \end{equation}

    By the triangle inequality and \Cref{assumption:relaxed-error}, the distance from the reference point to the current approximate iterate $x_k$ satisfies:
    \begin{equation}
        \norm{z_k - x_k} \le \norm{z_k - y_k} + \norm{y_k - x_k} \le r + \bar{\varepsilon} = t.
    \end{equation}
    Thus, $z_k \in B_t(x_k)$, meaning $z_k$ is a feasible point for the broximal subproblem solved at iteration $k$. 
    Because $y_{k+1}$ is defined as the exact minimizer of $\phi$ over this ball, we have:
    \begin{equation}
    \label{eq:general-y-and-z}
        \phi(y_{k+1}) \le \phi(z_k).
    \end{equation}

    Plugging \eqref{eq:general-y-and-z} into \eqref{eq:general-convexity} provides a one-step contraction bound for the sequence $\{y_k\}$:
    \begin{equation}
        \label{eq:general-one-step}
        \phi(y_{k+1}) - \phi_\star \le \left(1 - \frac{r}{\norm{y_k^\star - y_k}}\right) (\phi(y_k) - \phi_\star).
    \end{equation}

    To bound the denominator, we recall that, for every finite $1\le k\le K_t$, the point $y_k$ is generated by an exact broximal step from $x_{k-1}$.
    By \Cref{lem:exact-one-step}, the distance to the set of minimizers satisfies 
    \begin{equation}
        \dist(y_k, \cX_\star) \le \sqrt{D_{k-1}^2 - t^2} \le D_{k-1}.
    \end{equation}
    Furthermore, from the distance recursion in \Cref{lem:dist-recursion}(i), we know 
    \begin{equation}
        D_i \le D_{i-1} + \varepsilon_{i-1} \le D_{i-1} + \bar{\varepsilon} .
    \end{equation}
    By induction, it follows that $D_{k-1} \le D_0 + (k-1)\bar{\varepsilon}$ for every finite $1\le k\le K_t$. 
    Since $y_0 = x_0$ yields $\norm{y_0^\star - y_0} = D_0$, it follows that for every finite $0\le k\le K_t$, the distance to the set of minimizers is bounded by:
    \begin{equation}
    \label{eq:general-distance-bound}
        \norm{y_k^\star - y_k} = \dist(y_k, \cX_\star) \le D_0 + k \bar{\varepsilon}.
    \end{equation}

    Let $V_k \eqdef \phi(y_k) - \phi_\star$ denote the suboptimality gap. 
    Chaining the one-step contractions from \eqref{eq:general-one-step} for $k=0, \dots, N-1$ with $N \le K_t$, and applying the distance bound \eqref{eq:general-distance-bound} gives:
    \begin{equation}
        \label{eq:suboptimality-product}
        V_N \le V_0 \prod_{k=0}^{N-1} \left(1 - \frac{r}{D_0 + k \bar{\varepsilon}}\right).
    \end{equation}

    Applying the standard inequality $1-x \le \exp(-x)$ for $x \in (0,1)$ yields:
    \begin{equation}
        V_N \le V_0 \exp\left(-\sum_{k=0}^{N-1} \frac{r}{D_0 + k \bar{\varepsilon}}\right).
    \end{equation}
    Bounding the sum from below with its corresponding integral gives:
    \begin{equation}
        \sum_{k=0}^{N-1} \frac{r}{D_0 + k \bar{\varepsilon}} \ge \int_0^N \frac{r}{D_0 + x\bar{\varepsilon}} \, dx = \frac{r}{\bar{\varepsilon}} \ln\left(1 + \frac{N\bar{\varepsilon}}{D_0}\right).
    \end{equation}
    Substituting this into the exponential bound establishes the claim:
    \begin{equation}
        V_N \le V_0 \left(1 + \frac{N \bar{\varepsilon}}{D_0}\right)^{-\frac{r}{\bar{\varepsilon}}}.
    \end{equation}
\end{proof}

\begin{remark}[Convergence limit in the exact case]
Taking the limit of \eqref{eq:general-convex-rate} as $\bar{\varepsilon}\to0$ yields the bound $\phi(y_N)-\phi_\star\le\exp(-tN/D_0)(\phi(y_0)-\phi_\star)$.
For exact steps, setting $\bar{\varepsilon}=0$ and $r=t$ in \eqref{eq:suboptimality-product} gives the sharper bound
\[
    \phi(y_N)-\phi_\star
    \le\left(1-\frac{t}{D_0}\right)^N(\phi(y_0)-\phi_\star).
\]
This is the geometric decrease established in the proof of \citet[Theorem~22]{carmon2020acceleration} for the exact ball-oracle iteration with strictly convex objectives.
\citet[Theorem~3.1]{richtarik2026sharper} prove the same bound for \algname{Exact BPM} under \Cref{assumption:main}, which improves on the bound $(1+t/D_0)^{-N}$ of \citet[Corollary~8.2]{gruntkowska2025ball}.
\end{remark}

Writing $\bar{\varepsilon} = qt$ with $q \in (0,1)$, so that $r/\bar{\varepsilon} = (1-q)/q$, \Cref{theorem:two-step-analysis} gives the rate $\mathcal{O}(N^{-(1-q)/q})$; for example, $q = 1/4$, $1/2$, and $3/4$ give $\mathcal{O}(N^{-3})$, $\mathcal{O}(N^{-1})$, and $\mathcal{O}(N^{-1/3})$.
The exponent depends only on the relative error, and since $(1+N\bar{\varepsilon}/D_0)^{-r/\bar{\varepsilon}} \le (D_0/(N\bar{\varepsilon}))^{r/\bar{\varepsilon}}$, the constant hidden in $\mathcal{O}(\cdot)$ is $(D_0/\bar{\varepsilon})^{(1-q)/q}(\phi(y_0)-\phi_\star)$.
By contrast, the bounds of standard first-order methods scale with function-specific constants: subgradient methods for $G$-Lipschitz functions achieve $\mathcal{O}(GD_0/\sqrt{N})$ \citep{shor1985minimization, nesterov2018lectures}, gradient descent for $L$-smooth functions achieves $\mathcal{O}(LD_0^2/N)$ \citep{nesterov2018lectures}, and accelerated methods improve this to $\mathcal{O}(LD_0^2/N^2)$ \citep{nesterov1983method}.
Our bound requires no such regularity.
This feature is inherited from \algname{Exact BPM} \citep[Corollary~8.2]{gruntkowska2025ball} and survives persistent approximation errors.

\subsection{Bounded sublevel sets}

The sublinear rate of \Cref{theorem:two-step-analysis} can be sharpened to a linear rate if the sublevel sets of $\phi$ are bounded, which confines the reference points $y_k$ to a bounded region.

\begin{definition}[Bounded sublevel sets]
\label{def:bounded-sublevel}
    The function $\phi$ has bounded sublevel sets if for every $\alpha \in \R$, the corresponding sublevel set $\{y \in \R^d : \phi(y) \le \alpha\}$ is bounded. 
    We denote the sublevel set of the initial point $x_0$ by $L_{x_0}$ and its largest distance to the set of minimizers by $D_{\max}$:
    \begin{equation}
    \label{eq:l0-def}
        L_{x_0} \eqdef \{y \in \R^d : \phi(y) \le \phi(x_0)\},
    \end{equation}
    \begin{equation}
    \label{eq:dmax-def}
        D_{\max} \eqdef \sup_{y \in L_{x_0}} \dist(y, \cX_\star).
    \end{equation}
\end{definition}

\begin{restatable}{theorem}{sublevel}
\label{theorem:bounded-sublevel}
    Let $\{x_k\}_{k\ge 0}$ be the sequence generated by \Cref{alg:inexact-brox} and $\{y_k\}_{k\ge 0}$ the auxiliary sequence defined in \eqref{eq:y-seq}, and suppose \Cref{assumption:main,assumption:relaxed-error} hold. 
    Assume $D_0 > t$ and $\bar{\varepsilon} > 0$. 
    If $\phi$ has bounded sublevel sets (\Cref{def:bounded-sublevel}), then for any iteration $N \le K_t$, the suboptimality gap satisfies:
    \begin{equation}
    \label{eq:linear-rate-sublevel}
        \phi(y_N) - \phi_\star \le \min \left\{ \left(1 + \frac{N \bar{\varepsilon}}{D_0}\right)^{-\frac{r}{\bar{\varepsilon}}}, \; \left(1 - \frac{r}{D_{\max}}\right)^N \right\} (\phi(y_0) - \phi_\star),
    \end{equation}
    where $r$ is the inner radius defined in \eqref{eq:r-and-R}, and $D_{\max}$ is given by \eqref{eq:dmax-def}.
\end{restatable}

When $D_0 \ll D_{\max}$, the first bound in \eqref{eq:linear-rate-sublevel} can be much smaller in early iterations.
For large $N$, the second bound gives linear convergence with factor $1 - r/D_{\max}$.
This rate is no longer independent of the objective, since $D_{\max}$ depends on the sublevel set $L_{x_0}$ (cf. \Cref{tab:sec4}).
The theorem also shows that $K_t<\infty$, although it gives no bound on $K_t$.
Indeed, if $K_t=\infty$, then $\phi(y_N)\to\phi_\star$; since $L_{x_0}$ is compact and $\phi$ is lower semicontinuous, this forces $\dist(y_N,\cX_\star)\to0$, contradicting $\dist(y_N,\cX_\star)>r$ for all $N<K_t$ (see the proof of \Cref{theorem:two-step-analysis}).

\subsection{Strongly convex objectives}

Under strong convexity, the reference values $\phi(y_k)$ decrease by at least a fixed amount in every step before tube entry, which yields finite termination.
Our next result expresses the linear bound of \Cref{theorem:bounded-sublevel} in terms of the strong-convexity parameter and adds this absolute-decrease estimate.

\begin{definition}[Strong convexity]
\label{def:strong-convexity}
    Under \Cref{assumption:main}, the function $\phi$ is $\mu$-strongly convex with parameter $\mu > 0$ if for all $x, y \in \dom \phi$ and $\alpha \in [0,1]$, the following holds:
    \begin{equation}
    \label{eq:strong-convexity-def}
        \phi(\alpha x + (1-\alpha)y) \le \alpha \phi(x) + (1-\alpha)\phi(y) - \frac{\mu}{2}\alpha(1-\alpha)\norm{x - y}^2.
    \end{equation}
\end{definition}

\begin{restatable}{theorem}{strongconvex}
\label{theorem:strong-convexity-rate}
    Let $\{x_k\}_{k\ge 0}$ be the sequence generated by \Cref{alg:inexact-brox} and $\{y_k\}_{k\ge 0}$ the auxiliary sequence defined in \eqref{eq:y-seq}, and suppose \Cref{assumption:main,assumption:relaxed-error} hold. 
    Assume $D_0 > t$ and $\bar{\varepsilon} > 0$. 
    If $\phi$ is $\mu$-strongly convex, then for any iteration $N \le K_t$, the suboptimality gap satisfies:
    \begin{equation}
    \label{eq:strong-convex-rate}
        \phi(y_N) - \phi_\star \le \min \left\{ \left(1 + \frac{N \bar{\varepsilon}}{D_0}\right)^{-\frac{r}{\bar{\varepsilon}}} V_0, \; \left(1 - r\sqrt{\frac{\mu}{2V_0}}\right)^N V_0, \; V_0 - N \frac{\mu r^2}{2} \right\},
    \end{equation}
    where $V_0 = \phi(y_0) - \phi_\star$, and $r$ is the inner radius defined in \eqref{eq:r-and-R}.
\end{restatable}

Which of the three bounds is smallest depends on the parameters and on $N$.
Since the gap $\phi(y_N) - \phi_\star$ is nonnegative, the constant decrease $\mu r^2/2$ cannot continue indefinitely, which gives the following bound on the hitting time.

\begin{restatable}[Finite hitting time]{corollary}{strongconvexhit}
\label{cor:strong-convex-hitting}
    Under the conditions of \Cref{theorem:strong-convexity-rate}, the hitting time $K_t$ to the $t$-tube $\mathcal{T}_t(\cX_\star)$ is bounded by:
    \begin{equation}
    \label{eq:strong-convex-hitting-time}
        K_t \le \left\lceil \frac{2(\phi(y_0) - \phi_\star)}{\mu r^2} \right\rceil.
    \end{equation}
\end{restatable}

Thus, under strong convexity, the method still reaches the target neighborhood in finitely many steps when each subproblem is solved only up to an error below the ball radius.

\subsection{Discussion}
\label{sub:discussion}

\begin{table}[t]
\caption{\textbf{Suboptimality bounds.}
Here $D_0 = \dist(x_0, \cX_\star)$ is the initial distance to the set of minimizers, $V_0 = \phi(y_0) - \phi_\star = \phi(x_0) - \phi_\star$ is the initial suboptimality gap, and $r = t - \bar{\varepsilon}$ is the inner radius \eqref{eq:r-and-R}.
The results of this section hold under the relaxed error regime of \Cref{assumption:relaxed-error} for all iterations before entering the $t$-tube, $N \le K_t$.
For \Cref{theorem:bounded-sublevel,theorem:strong-convexity-rate}, the table shows only the bound that uses the additional assumption; the bound of \Cref{theorem:two-step-analysis} holds as well.}
\label{tab:sec4}
\small
\centering
\begin{tabular}{cccc}
\toprule
\textbf{Result} & \textbf{Assumptions} & \textbf{\makecell{Bound on \\ Suboptimality Gap}} & \textbf{\makecell{Independent of \\ Regularity}} \\
\midrule
\makecell{Theorem 3.1 \\ \citep{richtarik2026sharper}} & \makecell{Convexity, \\ exact steps} & $\left(1 - \frac{t}{D_0}\right)^{N} V_0$ & \checkmarkgreen \\
\addlinespace
\Cref{theorem:two-step-analysis} & \makecell{Convexity} & $\left(1 + \frac{N \bar{\varepsilon}}{D_0}\right)^{-\frac{r}{\bar{\varepsilon}}} V_0$ & \checkmarkgreen \\
\addlinespace
\Cref{theorem:bounded-sublevel} & \makecell{Convexity, \\ bounded sublevel sets} & $\left(1 - \frac{r}{D_{\max}}\right)^{N} V_0$ & \crossmarkred \\
\addlinespace
\Cref{theorem:strong-convexity-rate} & \makecell{$\mu$-strong convexity} & $V_0 - N \frac{\mu r^2}{2}$ & \crossmarkred \\
\bottomrule
\end{tabular}
\end{table}

We summarize the guarantees of this section in \Cref{tab:sec4}.
Under convexity alone, \Cref{theorem:two-step-analysis} gives a sublinear bound for the reference function values, valid up to tube entry and independent of the regularity of the objective.
In comparison, \algname{Exact BPM} admits a linear bound under convexity alone \citep{gruntkowska2025ball, richtarik2026sharper}.
With the additional assumption of bounded sublevel sets, we obtain a linear bound involving $D_{\max}$.
Under strong convexity, \Cref{cor:strong-convex-hitting} additionally supplies an explicit finite bound on $K_t$ and thus a complexity guarantee for reaching the $t$-tube.
This bound depends explicitly on the strong-convexity parameter $\mu$.

The guarantees of this section have two limitations.
First, they bound the suboptimality $\phi(y_N) - \phi_\star$ along the auxiliary sequence $\{y_k\}$ of \emph{exact} broximal points, rather than along the iterates $\{x_k\}$ actually produced by \Cref{alg:inexact-brox}. 
This distinction is unavoidable under our assumptions: the bound $\norm{x_k - y_k} \le \bar{\varepsilon}$ does not control $\phi(x_k) - \phi(y_k)$ for a general proper, closed, extended-valued convex $\phi$. 
Indeed, $x_k$ may lie arbitrarily close to $y_k$ yet satisfy $\phi(x_k) = +\infty$, so a suboptimality guarantee for the produced iterate would be false without additional regularity such as Lipschitz continuity or smoothness. 
The point $y_{k+1}$ is the exact target of the approximate subproblem solved at iteration $k$, and the two sequences remain linked by $\norm{x_{k+1} - y_{k+1}} \le \varepsilon_k$.

Second, these bounds hold only for $N \le K_t$, so they describe progress \emph{until} the trajectory enters the $t$-tube.
Except under strong convexity (\Cref{cor:strong-convex-hitting}), we have no a priori bound on $K_t$, so the rates cannot in general be converted into the number of iterations needed to reach $\phi(y_N) - \phi_\star \le \varepsilon$.

We next bound the suboptimality of the computed iterates after they enter the target neighborhood, in terms of the behavior of $\phi$ near the set of minimizers.
For $\rho \ge 0$, define the \emph{tube gap}
\begin{equation}
    \label{eq:tube-gap}
    \omega(\rho) \eqdef \sup\bigl\{\phi(x) : x \in \mathcal{T}_\rho(\cX_\star)\bigr\} - \phi_\star,
\end{equation}
the worst-case suboptimality of $\phi$ over the closed $\rho$-tube $\mathcal{T}_\rho(\cX_\star)$.
It is nondecreasing, with $\omega(0) = 0$.
The value $\omega(\rho)$ may be $+\infty$, in which case the bounds below are vacuous.

\begin{proposition}[Accuracy at tube entry]
\label{prop:tube-accuracy}
Let \Cref{assumption:main,assumption:relaxed-error} hold and suppose the hitting time is finite, $K_t < \infty$.
\begin{enumerate}[label=(\roman*)]
    \item The terminal iterate of \Cref{alg:inexact-brox} satisfies $\phi(x_{K_t}) - \phi_\star \le \omega(t)$.
    \item Consider the untruncated iteration, and assume that, for every $k\ge0$, there exists $y_{k+1}\in\broxSS{\phi}{t}{x_k}$ such that $\norm{x_{k+1}-y_{k+1}}\le\varepsilon_k$.
    Then $\dist(x_m,\cX_\star)\le\bar{\varepsilon}$, and hence $\phi(x_m)-\phi_\star\le\omega(\bar{\varepsilon})$, for every $m>K_t$.
\end{enumerate}
\end{proposition}

\begin{proof}
(i) By definition of $K_t$, the iterate satisfies $\dist(x_{K_t}, \cX_\star) \le t$, so $x_{K_t} \in \mathcal{T}_t(\cX_\star)$ and $\phi(x_{K_t}) - \phi_\star \le \omega(t)$ by \eqref{eq:tube-gap}. 
Truncation freezes the sequence, so the same bound holds for every $x_m$ with $m \ge K_t$.
(ii) Suppose $\dist(x_k, \cX_\star) \le t$. Then $B_t(x_k)$ contains a minimizer, so $\min_{y \in B_t(x_k)} \phi(y) = \phi_\star$ and the exact broximal point $y_{k+1}\in\broxSS{\phi}{t}{x_k}$ lies in $\cX_\star$. 
Hence $\dist(x_{k+1}, \cX_\star) \le \norm{x_{k+1} - y_{k+1}} \le \varepsilon_k \le \bar{\varepsilon} < t$, and the hypothesis holds again at $k+1$. 
By induction, $\dist(x_m, \cX_\star) \le \bar{\varepsilon}$ for all $m > K_t$, and $\phi(x_m) - \phi_\star \le \omega(\bar{\varepsilon})$ by \eqref{eq:tube-gap}.
\end{proof}

\section{Numerical illustration}
\label{sec:num-exp}

We illustrate the convergence behavior of \Cref{alg:inexact-brox} on a simple quadratic objective in $\R^2$,
\begin{equation*}
   \phi(x) = x_1^2 +3x_2^2, \quad x_\star = (0, 0)^\top,
\end{equation*}
starting from $x_0 = (10, 10)^\top$ with radius $t=1$. 

To evaluate the theoretical guarantees, we consider an adversarial approximation operator that, within the $\bar{\varepsilon}$-ball centered at $y_{k+1}$, selects the point maximizing the distance to the optimum:
\begin{equation*}
    x_{k+1} = y_{k+1} + \bar{\varepsilon} \frac{y_{k+1}}{\norm{y_{k+1}}}.
\end{equation*}
The adversary thus pushes the iterates as far away from the optimum as possible at each step.

\begin{figure}[htbp]
  \centering
  \includegraphics[width=\textwidth]{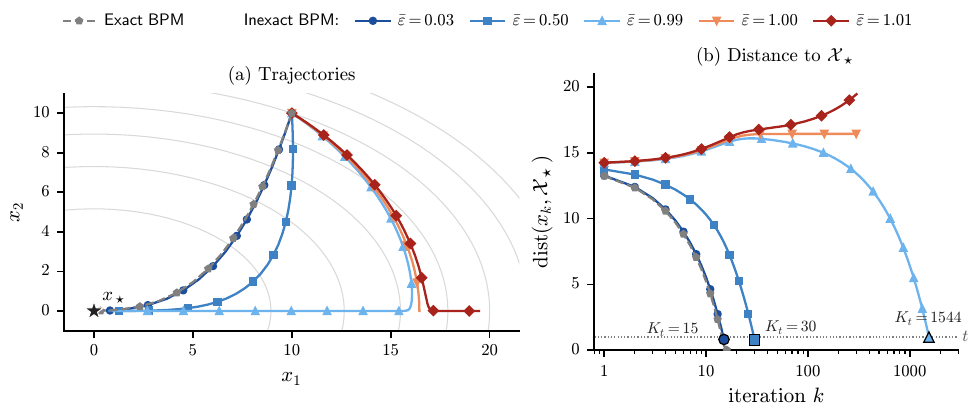}
  \caption{\textbf{Robustness of \algname{Inexact BPM} under adversarial perturbations.} \algname{Exact BPM} and \algname{Inexact BPM} with radius $t=1$ on $\phi(x) = x_1^2 +3x_2^2$ from $x_0=(10,10)^\top$, using an adversarial operator that maximizes the distance to the optimum within $\bar{\varepsilon} \in \{0.03, 0.50, 0.99, 1.00, 1.01\}$. (a)~Trajectories. (b)~Distance to the optimum against the iteration count (log scale); black-edged markers indicate the hitting time $K_t$. For every tested $\bar{\varepsilon} < t$, the iterates reach the $t$-tube although the adversary pushes them outward. The runs with $\bar{\varepsilon} \ge t$ never enter the tube and are stopped after $300$ iterations.}
  \label{fig:adversarial}
\end{figure}

\Cref{fig:adversarial} shows the trajectories of \algname{Exact BPM} and \algname{Inexact BPM} for $\bar{\varepsilon} \in \{0.03, 0.50, 0.99, 1.00, 1.01\}$, together with the distance to the optimum at each iteration.
\algname{Exact BPM} reaches $x_\star$ after $16$ steps.
For $\bar{\varepsilon}<t=1$, \algname{Inexact BPM} reaches the $t$-tube, as guaranteed by \Cref{cor:strong-convex-hitting}, since the quadratic objective is strongly convex with $\mu=2$.
The hitting times are $K_t = 15$, $30$, and $1544$ for $\bar{\varepsilon} = 0.03$, $0.50$, and $0.99$, well below the corresponding bounds $426$, $1600$, and $4\cdot10^6$ of \Cref{cor:strong-convex-hitting}.
For $\bar{\varepsilon} = 0.03$, the condition $\bar{p} > 0$ of \Cref{theorem:uniform-small-errors} also holds, since $D_0 - \sqrt{D_0^2 - t^2} \approx 0.035$ for $D_0 = 10\sqrt{2}$, and that theorem gives $K_t \le 1310$.
Convergence slows as $\bar{\varepsilon}$ approaches $t$, in line with the bounds, which degrade as the inner radius $r=t-\bar{\varepsilon}$ tends to zero.
For $\bar{\varepsilon} \ge t$, \Cref{assumption:relaxed-error} fails.
With $\bar{\varepsilon} = t$, the trajectory is first pushed outward and then stalls; with $\bar{\varepsilon} > t$, the adversary pushes it outward indefinitely.

\section{Conclusion}
\label{sec:conclusion}

We introduced \algname{Inexact BPM} and gave, to our knowledge, the first convergence analysis of the Broximal Point Method with inexact broximal steps.
We measure inexactness by the distance to the exact broximal point and make no assumption on how the approximation is obtained.
Under summable or uniformly small errors, we bounded the hitting time of the $t$-tube (\Cref{theorem:summable-progress-guarantee,theorem:uniform-small-errors}).
Under the weaker assumption that the errors are uniformly bounded by a constant below the ball radius ($\bar{\varepsilon} < t$), we proved sublinear convergence in function value for general convex objectives (\Cref{theorem:two-step-analysis}), which sharpens to a linear rate under bounded sublevel sets (\Cref{theorem:bounded-sublevel}) and yields an explicit bound on the hitting time under strong convexity (\Cref{cor:strong-convex-hitting}).
In the general convex setting, as for \algname{Exact BPM}, these bounds do not depend on the condition number or on any other regularity property of the objective.
In the general convex and bounded-sublevel settings, however, they are function-value guarantees for the exact broximal points that hold until the $t$-tube is reached, not complexity bounds for reaching $\varepsilon$-optimality with the computed iterates.

Since our analysis imposes no structure on how the approximations are computed, it applies to any inner solver whose output satisfies the error conditions.
Designing such solvers, and certifying their accuracy in practice, is left to future work.

\subsubsection*{Acknowledgments}
The research reported in this publication was supported by funding from King Abdullah University of Science and Technology (KAUST): i) KAUST Baseline Research Scheme, ii) CRG Grant ORFS-CRG12-2024-6460, and iii) Center of Excellence for Generative AI, under award number 5940.

\medskip

\bibliographystyle{abbrvnat}
\bibliography{references}

\newpage
\appendix

\section{Proofs}
\label{sec:proofs}

\subsection{\proofof{Lemma}{lem:exact-one-step}}

\broxlemma*

\begin{proof}
We adapt the argument of \citet[Theorem~D.1]{gruntkowska2025ball}.
Since $\phi$ is lower semicontinuous and $B_t(x)$ is compact, the minimum is attained; it is finite because $B_t(x)\cap\dom\phi\neq\emptyset$ and $\phi\ge\phi_\star$.
An interior solution would be a local minimizer of $\phi$ and, by convexity, a global minimizer, contradicting $B_t(x)\cap\cX_\star=\emptyset$.
Hence every solution lies on the boundary of the ball.
If two distinct solutions existed, their midpoint would also be a solution by convexity, but would lie in the interior of the Euclidean ball.
Thus the solution $y$ is unique and satisfies $\norm{y-x}=t$.

Fix $x_\star\in\cX_\star$. Since $\phi(y)>\phi_\star$, convexity implies
\[
    \phi\bigl(y+\lambda(x_\star-y)\bigr)<\phi(y),\qquad 0<\lambda\le1.
\]
Optimality of $y$ therefore gives $\norm{y+\lambda(x_\star-y)-x}^2>t^2$.
Expanding, dividing by $\lambda$, and letting $\lambda\downarrow0$ yields $\langle x-y,y-x_\star\rangle\ge0$. Hence
\[
    \norm{x-x_\star}^2=t^2+\norm{y-x_\star}^2+2\langle x-y,y-x_\star\rangle \ge t^2+\norm{y-x_\star}^2.
\]
\end{proof}

\subsection{\proofof{Lemma}{lem:dist-recursion}}
\label{sub:proof_of_lem:dist-recursion}

\recursionlemma*

\begin{proof}
    Assume first that $D_k > t$, and let
    \begin{equation*}
        y_{k+1} = \broxSS{\phi}{t}{x_k}.
    \end{equation*}
    By the one-step contraction property (\Cref{lem:exact-one-step}), we have 
    \begin{equation*}
        \norm{y_{k+1} - z}^2 \le \norm{x_k - z}^2 - t^2
    \end{equation*}
    for every minimizer $z \in \cX_\star$. 
    Let $x_k^\star \in \cX_\star$ be the projection of $x_k$ onto $\cX_\star$, so that $\norm{x_k - x_k^\star} = \dist(x_k, \cX_\star)$. Substituting $z = x_k^\star$ yields:
    \begin{equation*}
        \dist(y_{k+1}, \cX_\star)^2 \le \norm{y_{k+1} - x_k^\star}^2 \le \norm{x_k - x_k^\star}^2 - t^2 = D_k^2 - t^2,
    \end{equation*}
    and therefore,
    \begin{equation}
    \label{eq:y_dist_bound}
        \dist(y_{k+1}, \cX_\star) \le \sqrt{D_k^2 - t^2}.
    \end{equation}
    
    By \eqref{eq:approximation-error}, the approximate iterate satisfies $\norm{x_{k+1} - y_{k+1}} \le \varepsilon_k$.
    The triangle inequality now yields
    \begin{equation*}
        D_{k+1} = \dist(x_{k+1}, \cX_\star) \le \dist(y_{k+1}, \cX_\star) + \norm{x_{k+1} - y_{k+1}}.
    \end{equation*}
    Applying \eqref{eq:y_dist_bound} and \eqref{eq:approximation-error} to the right-hand side gives:
    \begin{equation*}
        D_{k+1} \le \sqrt{D_k^2 - t^2} + \varepsilon_k,
    \end{equation*}
    which proves part (i).

    If $D_k \le t$, the algorithm sets $x_{k+1} = x_k$. 
    Thus, $D_{k+1} = D_k \le t$, establishing part (ii).
\end{proof}

\subsection{\proofof{Theorem}{theorem:summable-progress-guarantee}}
\label{sub:proof_of_theorem:summable-progress-guarantee}

\summable*

\begin{proof}
    If $D_0 \le t$, then $K_t = 0$ and the claimed bound holds trivially. 

    Assume henceforth that $D_0 > t$.
    By definition, $D_k > t$ for all $k < K_t$, so \Cref{lem:dist-recursion}(i) applies at every such step:
    \begin{equation}
    \label{eq:one-step-applied-sum}
        D_{k+1}^2 \le D_k^2 - t^2 + \varepsilon_k^2 + 2\varepsilon_k\sqrt{D_k^2 - t^2}, \qquad \forall k < K_t.
    \end{equation}
    We upper bound the square root term by $D_k$, so
    \begin{equation}
    \label{eq:summable-eq1}
        D_{k+1}^2 \le D_k^2 - t^2 + \varepsilon_k^2 + 2\varepsilon_k D_k, \qquad \forall k < K_t.
    \end{equation}
    By \Cref{lem:dist-recursion}(i), $D_{i+1}\le D_i+\varepsilon_i$ for every $i<K_t$.
    Hence $D_k\le D_0+\sum_{i=0}^{k-1}\varepsilon_i$ for every finite $k\le K_t$.
    Substituting this bound into \eqref{eq:summable-eq1} yields
    \[
        D_{k+1}^2
        \le D_k^2-t^2+\varepsilon_k^2
        +2\varepsilon_k\left(D_0+\sum_{i=0}^{k-1}\varepsilon_i\right).
    \]
    For any finite integer $N\le K_t$, summing over $k=0,\ldots,N-1$ and using
    \[
        \sum_{k=0}^{N-1}\varepsilon_k^2
        +2\sum_{k=0}^{N-1}\varepsilon_k\sum_{i=0}^{k-1}\varepsilon_i
        =\left(\sum_{k=0}^{N-1}\varepsilon_k\right)^2
    \]
    gives
    \[
        D_N^2
        \le\left(D_0+\sum_{k=0}^{N-1}\varepsilon_k\right)^2-Nt^2
        \le(D_0+\Sigma)^2-Nt^2.
    \]
    Since $D_N^2\ge0$, we obtain $Nt^2\le(D_0+\Sigma)^2$.
    If $K_t=\infty$, choosing a finite integer $N>(D_0+\Sigma)^2/t^2$ gives a contradiction.
    Therefore, $K_t<\infty$, and taking $N=K_t$ establishes the claimed bound.
\end{proof}

\subsection{\proofof{Theorem}{theorem:uniform-small-errors}}
\label{sub:proof_of_theorem:uniform-small-errors}

\uniformsmall*

\begin{proof}

    We first verify that every active subproblem is well-defined.
    The condition $\bar p>0$ implies $\bar\varepsilon<D_0-\sqrt{D_0^2-t^2}<t$.
    The initial ball intersects $\dom\phi$ because $x_0\in\dom\phi$.
    Suppose that the subproblem at iteration $k<K_t$ is well-defined, and let $y_{k+1}$ denote its exact solution.
    Then $y_{k+1}\in\dom\phi$ and $\norm{x_{k+1}-y_{k+1}}\le\varepsilon_k\le\bar\varepsilon<t$, so $B_t(x_{k+1})\cap\dom\phi\neq\emptyset$.
    Thus every active subproblem is well-defined by induction.
    Since $D_k>t$ for all $k<K_t$, \Cref{lem:dist-recursion}(i) yields:
    \begin{equation}
    \label{eq:one-step-applied}
        D_{k+1}^2 \le D_k^2 - t^2 + \varepsilon_k^2 + 2\varepsilon_k\sqrt{D_k^2 - t^2}, \qquad \forall k < K_t.
    \end{equation}
    Since $\varepsilon_k \le \bar{\varepsilon}$ and the right-hand side of \eqref{eq:one-step-applied} is monotonically increasing with respect to $\varepsilon_k \ge 0$, we can substitute the uniform error bound:
    \begin{equation}
    \label{eq:one-step-applied-uniform}
        D_{k+1}^2 \le D_k^2 - t^2 + \bar{\varepsilon}^2 + 2\bar{\varepsilon}\sqrt{D_k^2 - t^2}, \qquad \forall k < K_t.
    \end{equation}
    Evaluating \eqref{eq:one-step-applied-uniform} at $k=0$ under the assumption that $\bar{p} > 0$ yields $D_1^2 \le D_0^2 - \bar{p} < D_0^2$. 
    If $D_k\le D_0$ and $k<K_t$, then \eqref{eq:one-step-applied-uniform} gives $D_{k+1}^2\le D_k^2-\bar p<D_k^2$.
    Induction therefore shows that the sequence $\{D_k\}_{k=0}^{K_t}$ is decreasing and that $D_k\le D_0$.
    
    Substituting the upper bound $D_k \le D_0$ into the positive square root term of \eqref{eq:one-step-applied-uniform} yields a constant minimum progress per step:
    \begin{equation*}
        D_{k+1}^2 \le D_k^2 - \left(t^2 - \bar{\varepsilon}^2 - 2\bar{\varepsilon}\sqrt{D_0^2 - t^2}\right) = D_k^2 - \bar{p}.
    \end{equation*}
    Fix any finite integer $N\le K_t$. 
    Summing the preceding inequality from $k=0$ to $N-1$ gives
    \[
        N\bar p\le D_0^2.
    \]
    If $K_t=+\infty$, we could choose a finite integer $N>D_0^2/\bar p$, obtaining a contradiction. 
    Hence $K_t<+\infty$, and taking $N=K_t$ gives the claimed bound.
\end{proof}

\subsection{\proofof{Theorem}{theorem:bounded-sublevel}}
\label{sub:proof_of_theorem:bounded-sublevel}

\sublevel*
\begin{proof}
    Since $L_{x_0}$ is bounded, $D_{\max} = \sup_{y \in L_{x_0}} \dist(y, \cX_\star) < \infty$. 
    For every active iteration $k<K_t$, combining \eqref{eq:general-y-and-z} and \eqref{eq:general-convexity} gives
    \[
        \phi(y_{k+1}) \le \phi(z_k) \le \phi(y_k).
    \]
    Hence, by induction, $\phi(y_k)\le \phi(y_0)=\phi(x_0)$ for every $k\le K_t$.
    Therefore, $y_k\in L_{x_0}$ for every $k\le K_t$, and consequently
    \[
        \norm{y_k-y_k^\star} = \dist(y_k,\cX_\star) \le D_{\max}.
    \]

    We now have two independent upper bounds on the distance to the set of minimizers for all $k \le K_t$: 
    the global bound $\norm{y_k - y_k^\star} \le D_{\max}$, and the accumulated error bound $\norm{y_k - y_k^\star} \le D_0 + k \bar{\varepsilon}$ derived in the proof of \Cref{theorem:two-step-analysis}.

    Substituting the accumulated error bound into the one-step contraction \eqref{eq:general-one-step} and chaining over $N$ steps recovers the sublinear rate from \Cref{theorem:two-step-analysis}:
    \begin{equation*}
        V_N \le V_0 \prod_{k=0}^{N-1} \left(1 - \frac{r}{D_0 + k\bar{\varepsilon}}\right) \le V_0 \left(1 + \frac{N \bar{\varepsilon}}{D_0}\right)^{-\frac{r}{\bar{\varepsilon}}}.
    \end{equation*}

    Alternatively, substituting the global bound $\norm{y_k - y_k^\star} \le D_{\max}$ into \eqref{eq:general-one-step} and chaining over $N$ steps yields the linear rate:
    \begin{equation*}
        V_N \le V_0 \prod_{k=0}^{N-1} \left(1 - \frac{r}{D_{\max}}\right) = V_0 \left(1 - \frac{r}{D_{\max}}\right)^N.
    \end{equation*}

    Since both upper bounds hold, the suboptimality gap $V_N$ is bounded by their minimum, completing the proof.
\end{proof}

\subsection{\proofof{Theorem}{theorem:strong-convexity-rate}}
\label{sub:proof_of_theorem:strong-convexity-rate}

\strongconvex*

\begin{proof}
    Because $\phi$ is $\mu$-strongly convex, its set of minimizers is a singleton, $\cX_\star = \{x_\star\}$, and its sublevel sets are bounded.
    Specifically, for any $y \in L_{x_0}$, strong convexity implies:
    \begin{equation*}
        \frac{\mu}{2}\norm{y - x_\star}^2 \le \phi(y) - \phi_\star \le \phi(y_0) - \phi_\star = V_0.
    \end{equation*}
    Thus, the largest distance to $x_\star$ within the initial sublevel set is bounded by $D_{\max} \le \sqrt{2V_0/\mu}$. 
    Substituting this $D_{\max}$ into the linear bound from \Cref{theorem:bounded-sublevel} yields the second term in the minimum. 
    The first term follows directly from \Cref{theorem:two-step-analysis}.

    As shown in the proof of \Cref{theorem:two-step-analysis}, $\norm{x_k - y_k} \le \bar{\varepsilon}$ implies $B_r(y_k) \subseteq B_t(x_k)$, and hence:
    \begin{equation*}
        \phi(y_{k+1}) \le \phi(w_k) \eqdef \min_{z \in B_r(y_k)} \phi(z),
    \end{equation*}
    where $w_k$ is a minimizer of $\phi$ over $B_r(y_k)$. 
    Since $y_k\in\dom\phi$, $r>0$, and $w_k\in\broxSS{\phi}{r}{y_k}$, all hypotheses of \citet[Theorem~D.2]{gruntkowska2025ball} are satisfied.
    Consequently, there exists $c_r(y_k)\ge0$ such that
    \[
        c_r(y_k)(y_k-w_k)\in\partial\phi(w_k).
    \]
    Applying $\mu$-strong convexity with this subgradient yields
    \begin{align*}
        \phi(y_k)
        &\ge \phi(w_k) + \left\langle c_r(y_k)(y_k-w_k), y_k-w_k \right\rangle + \frac{\mu}{2}\norm{y_k-w_k}^2\\
        &= \phi(w_k) + c_r(y_k)\norm{y_k-w_k}^2 + \frac{\mu}{2}\norm{y_k-w_k}^2\\
        &\ge \phi(w_k) + \frac{\mu}{2}\norm{y_k-w_k}^2.
    \end{align*}
    For every active iteration $k<K_t$,
    \[
        \dist(y_k,\cX_\star) \ge \dist(x_k,\cX_\star)-\norm{x_k-y_k} > t-\bar{\varepsilon} = r.
    \]
    Hence, $\cX_\star\cap B_r(y_k)=\emptyset$.
    By \citet[Theorem~D.1(ii)]{gruntkowska2025ball}, this boundary condition forces $\norm{y_k - w_k} = t - \bar{\varepsilon} = r$. 
    Combining these inequalities, we obtain a guaranteed absolute function value decrease of:
    \begin{equation*}
        \phi(y_{k+1}) \le \phi(w_k) \le \phi(y_k) - \frac{\mu}{2}r^2.
    \end{equation*}
    Applying this recursively over $N$ steps gives $\phi(y_N) - \phi_\star \le V_0 - N \frac{\mu r^2}{2}$. 
    Since all three bounds hold, the suboptimality gap is bounded by their minimum.
\end{proof}

\subsection{\proofof{Corollary}{cor:strong-convex-hitting}}
\label{sub:proof_of_cor:strong-convex-hitting}

\strongconvexhit*

\begin{proof}
    By \Cref{theorem:strong-convexity-rate}, for every finite integer $N\le K_t$,
    \[
        0 \le \phi(y_N)-\phi_\star \le V_0-N\frac{\mu r^2}{2}.
    \]
    Consequently,
    \[
        N\le \frac{2V_0}{\mu r^2} \qquad \text{for every finite }N\le K_t.
    \]
    If $K_t=+\infty$, we could choose a finite integer $N>2V_0/(\mu r^2)$, which would contradict the preceding inequality. 
    Hence, $K_t<+\infty$. 
    Taking $N=K_t$ now gives
    \[
        K_t\le \frac{2V_0}{\mu r^2} \le \left\lceil\frac{2V_0}{\mu r^2}\right\rceil,
    \]
    which proves the claim.
\end{proof}

\end{document}